\documentstyle[amssymb,amscd,10pt]{amsart}

\theoremstyle{definition}
\newtheorem{defn}{Definition}[section]
\theoremstyle{definition}
\newtheorem{thm}[defn]{Theorem}

\newtheorem{lem}[defn]{Lemma}
\theoremstyle{remark}
\newtheorem{rem}{Remark}
\theoremstyle{example}
\newtheorem{exam}{Example}

\begin{document}


\renewcommand\sp{\operatorname{Spec}}       
\newcommand\aut{\operatorname{Aut}}
\newcommand\g {{\operatorname{G_\lambda}}(O_D^n)}
\renewcommand\hom{\operatorname{Hom}}
\newcommand\gi {{\operatorname{G_{\lambda+i}}}(O_D^n)}
\renewcommand\o{{{\cal O}}}      
 \newcommand\Z{{\Bbb Z}}    
\newcommand\f{{\Bbb F_q}}
\newcommand\s{\sigma}
\newcommand\w{\widehat}

\newcommand\limi{\varinjlim}

\newcommand\limp{\varprojlim}


\newcommand\iso{\East{\sim}{}}
\newcommand\co{\underset{A\sim{V^+}}\to{\underset {A}\to{\varprojlim}}}
\newcommand\M{\cal M}

\newcommand\punt {{\raise.4ex\hbox{\pun ^^O}}}
\newcommand\fum{
		\setbox4=\hbox{$\widetilde{\phantom{\widehat{\cal M}}}$}\ht4=15truept\dp4=0pt
		\setbox5=\hbox{${\hbox{$\widehat{\cal M}$}\atop\box4}$}\dp5=-5truept
		\raise-7truept\box5}


\newcommand\fu{\underline}

\title [ Drinfeld moduli schemes ]{  Drinfeld moduli schemes
  \\ and infinite Grassmannians}

\author[A. \'Alvarez]{A. \'Alvarez${}^*$}

\address{Departamento de Matem\'atica Pura y Aplicada \\
Universidad de Salamanca \\ Plaza de la Merced 1-4 \\
Salamanca 37008. Spain.}

\thanks{This work is partially supported by the CIGYT research
contract n. PB91-0188 \\$*$  Departamento de Matem\'atica Pura y Aplicada. 
Universidad de Salamanca}

\date{June 1997}

\maketitle 


\tableofcontents


\section{Introduction}

The aim of this paper is to construct an immersion of the  Drinfeld moduli schemes in a finite
product of infinite Grassmannians, such that they will be locally closed subschemes
of these Grassmannians which represent a kind of flag varieties. This construction is derived from
two   results: the first   is  that the moduli functor of vector bundles with an $\infty
'$-formal level structure (defined below) over a curve
$X$ is representable, in a natural way,  by a closed subscheme of the infinite Grassmannian.  
The second   is  an equivalence (see {\cite{D1}, {\cite{Mu}}, {\cite{BlSt}}, {\cite{LRSt}},
{\cite{C}}) between  Drinfeld
$A$-modules of rank
$n$ and elliptic sheaves, extended for level structures in {\cite{An}}.    

Let us detail these results. Let $X$ be a smooth, proper and geometrically irreducible curve over
a finite field
${{\Bbb F}_q}$, $\infty$ a rational point of $X$, $A=H^0(X-\infty,{{\cal O}_X})$, $D$ an
effective divisor over $Spec(A)$, and $S$ an arbitrary scheme over ${{\Bbb F}_q}$.  $\w {\Bbb
A}^\infty$ denotes the ring of adeles outside $\infty$ and 
$$\w {O}^\infty=\underset {I{\text{ ideal of }}A}   {\underset {I}   \varprojlim}A/I$$

  The moduli functor of vector
bundles of rank $n$, with a $\infty '$-formal level structure is:   
$$\w{\fu{ \cal M}}_\infty^n(S)=\underset
{\infty \notin supp(D)}   {\underset {D>0}
 \varprojlim}\fu {\cal M}^n_D(S)$$
 where

$$ \fu {{\cal M}}^n_D(S)= \left\{\vcenter{\hsize=7cm 
\parindent=0pt
$D$-level structures, $(M,f_D)$, over a locally free 
\par sheaf $M$ of rank $n$ over $X\times S$, up to isomorphism}\right\}$$

In this paper, we consider the moduli functor  \cite{D2}, {\cite{L}},  \cite{G}: 
$$\w{\fu{ \cal D }}_\infty^n =\underset {\infty \notin supp(D) } {\underset {D>0}
 \varprojlim}  \fu{ \cal D }^n_D $$
where $  \fu {\cal D}^n_D$ is the moduli functor of   Drinfeld $A$-modules endowed with
a $D$-level structure   with characteristic away from $supp(D)$. Or equivalently  elliptic
sheaves with a twisted $\infty '$-level structure (see Remark \ref{dri}).

The results will be obtained bearing in mind that the infinite Grassmannian scheme 
${\text{Gr}}((\w {\Bbb A}^\infty)^n,(\w {O}^\infty)^n)$ (\cite{BS}, \cite{ADKP}, \cite{AMP},
\cite{Gi})  represents the functor 

$\fu{\text{Gr}}((\w {\Bbb A}^\infty)^n,(\w {O}^\infty)^n)$ 
and that the Krichever map (\cite{M}, \cite{SW},\cite{Q}) given in  (\ref{kri})

$$\varphi: \w{\fu{ \cal M }}_\infty^n 
\to 
\fu {\text{Gr}}((\w {\Bbb A}^\infty)^n,(\w {O}^\infty)^n)$$

makes   $\w{\fu{ \cal M }}_\infty^n$ a closed subfunctor of
$\fu{\text{Gr}}((\w {\Bbb A}^\infty)^n,(\w {O}^\infty)^n)$. Therefore, it is representable by a
closed subscheme of
$Gr((\w {\Bbb A}^\infty)^n,(\w {O}^\infty)^n)$. In the literature  this morphism   is
considered when $S$ is a field, recently Ines Quandt extended it for noetherian schemes. In this
paper it is     considered for an arbitrary scheme, although fixing the curve $X$, (when
$S$  is the base field, those results are the Weil's uniformation of vector bundles,  in
\cite{BL} are settled   similar results but in another setting of infinite Grassmannian).

Putting together the Krichever
morphism with the equivalence between Drinfeld $A$-modules and elliptic sheaves,    we
obtain  in a natural way, a locally closed immersion of functors

$$ \psi:\w{\fu{ \cal D }}_\infty^n   \hookrightarrow
\underset {i=0}
 {\overset   { n-1} 
\prod} 
 \fu{\text{Gr}} ({\Omega^1}_{A/k}\underset {A}\otimes (\w {\Bbb
A}^\infty)^n,{\Omega^1}_{A/k}\underset {A}\otimes (\w {O}^\infty)^n)$$
 (${\Omega^1}_{A/k}$ is the sheaf of differentials over $A$)

 This infinite Grassmannian is a
slight modification   of  $\fu{\text{Gr}}((\w {\Bbb A}^\infty)^n,(\w {O}^\infty)^n)$.
Therefore, $ \w{\fu{ \cal D }}_\infty^n$ is representable by a locally closed
subscheme of a finite product of infinite Grassmannians, moreover it can been seen as a point of
an infinite flag variety, in this setting, one can see an analogy with Deligne-Lusztig variety for
the Coxeter element \cite{DL}.

Similar results could be obtained when 
$\underset {m\in {\Bbb N}}   
  \varprojlim  \fu{ \cal D }^n_{mD}$ 
and  
$\underset {m\in {\Bbb N}}   
 \varprojlim \fu {\cal M}^n_{mD}(S)$ 
are considered    instead of the projective limits:
$$\underset {\infty \notin supp(D)  } {\underset {D>0}
 \varprojlim}  \fu{ \cal D }^n_D\qquad\qquad \underset
{\infty \notin supp(D)  } {\underset {D>0}
 \varprojlim}\fu {\cal M}^n_{D}(S)$$  

  Let us now  briefly state the   contents of the different sections of this paper.
In the  second   and third  sections  we study the infinite Grassmannian and the Krichever
morphism over an arbitrary scheme  and conclude that
$\w{\fu{ \cal M }}_\infty^n$ is representable by  a
closed subscheme of  $ \fu{\text{Gr}} ((\w {\Bbb A}^\infty)^n,(\w {O}^\infty)^n)$. Using this, in
the fourth section, we show that the moduli functors of Drinfeld $A$-modules with level
structures are representable by locally closed schemes of a finite product of
infinite Grassmannians.

 
\section{Preliminaries on the Sato's infinite Grassmannian. \label{grass}}

Let $V$ be a vector space over a field $k$ and $A ,B$   vector subspaces of $V$.

\begin{defn}  
Two vector subspaces  $A ,B$ of
$V$   are said to be commensurable   if ${A+B}/{A\cap B}$ is a vector
space over $k$ of finite dimension.  \rm{\cite{T}}.
\end{defn}

Let us consider a collection, $F=\{{V_i}\}_{i \in I}$, of subspaces of $V$, such that $\underset
{i}
  \cap {V_i} ={0}$ and $V_i$ is commensurable to $V_j$, for all $i,j$,
with $F$ it is possible to define a Haussdorff topology on $V$ (Tate's topology): $F$ is a basis
of neighborhoods of $0$.

$A\sim F$ will denote an open subspace of $V$ commensurable with any  ${\{V_i}\}_{i \in I}$ (or
equivalently with some $V_i$).

\begin{defn}
 The completion of $V$ with respect to the
$F$-topology is defined by:
$$ \w V= \underset {A\sim F}   {\varprojlim}   V/A $$
\end{defn}

Analogously, given a vector subspace $B\subseteq V$ we can define the
completions of $B$ and $ V/B$ with respect to $B\cap F=\{{A\cap V_i}\}_{i \in I}$ and 
${B+F}/{B}=\{{B+V_i}/{B}\}_{i \in I}$,  respectively.

\begin{exam}
\begin{itemize}

\item $(V, F=\{0\})$; $V$ is complete.
\item $V=k((t)),$  $F=\{t^nk[[t]]\}_{n\in \Bbb Z}$; $V$ is complete.
\item Let $(X, \o_X)$ be a smooth, proper and
irreducible curve over the field $k$, and let $V$ be the ring of
adeles of the curve and $F=\{I\underset {p } {\prod} \w {O_p}\}_{I\subset \o_X}$
( $\w {O_p}$ being the ${ m}_p$-adic completion of the local
ring of $X$ in the point $p$) and $I$ ideals of $\o_X$; $V$ is complete with respect the
$F$-topology.
\end{itemize}
\end{exam}

\begin{defn}
Given a $k$-scheme $S$ and a vector
subspace $B\subseteq V$, we define the sheaves over $S$:
\begin{itemize}

\item  $\w V_S=
 \underset {A\sim F}   {\varprojlim} ( V / A\underset {k} {\otimes}{{\cal O}}_S).$

\item   $\widehat B_S=\underset {A\sim F}   {\varprojlim} (( {B} /{  A\cap B}) \underset {k }
{\otimes}{{\cal O}}_S)$.

\item $\widehat{( V/B)}_S=\underset {A\sim F}   {\varprojlim} (( {V}/{A+B})\underset {k}
 {\otimes}{{\cal O}}_S)$.
\end{itemize}
\end{defn}
We have the exact sequence (\cite{Ha},II.9.1 and II.9.2.1):

$$0\to \w B_S \to\w V_S \to\widehat{( V/B)}_S$$
 Remark that $\w V_S$ and  $\w B_S$ are not quasicoherent $\o_S$-modules but $\widehat{(
V/B)}_S$ it is. 

 A discrete vector subspace of $V$ is a
vector subspace, $  L\subseteq V$, such that $  L\cap V_i$ and
${V}/{  L+V_i}$ are $k$-vector spaces of finite dimension (for some $i \in I$).

\begin{defn}
 Given a $k$-scheme $S$, a discrete
submodule of $\w V_S $ is a subfunctor ${{\cal L}} \subset \w V_S $ in the category of
$S$-schemes, given by  a   quasicoherent sheaf of
$\o_S$-modules, verifying: for each
$s\in S$  
    there exists an open neighborhood
$U_s$ of
$s$ and an open commensurable $k$-vector subspace $B \sim F$ such that: ${\cal L}_{U_s}\cap \w
B_{U_s} $ is free of  type finite   and
$ {\w V_{U_s}} /{{\cal L}_{U_s}+ \w B_{U_s}}=0.$
\end{defn}

We define the Grassmannians functor over the category of $k$-schemes by:
 
$$\fu{\text{Gr}}(V,F)(S)= \left\{\vcenter{\hsize=6cm 
\parindent=0pt
discrete sub-$\o_S$-modules of
$\widehat{ V}_S $ with
\par   respect to the Tate's topology for $F$}\right\}$$

\begin{rem}   \rm{ Note that if $V$ is a finite dimensional
$k$-vector space and
$F=\{0\}$, then $\fu{\text{Gr}}(V,\{0\})$ is the usual Grassmannian
functor defined by Grothendieck \cite{{EGA}}. ( ${\cal L}$ is a subfuntor of $V\otimes \o_S$ in
the category of $S$-schemes if and only if $V\otimes {\o}_S/{{\cal L}}$ is flat.)}\end{rem}

\begin{thm} \rm{(\cite{ADKP},\cite{BS}, \cite{AMP}, \cite{Gi})}  
The functor $\fu{\text{Gr}}(V,F)$ is representable by a $k$-scheme
$\text{Gr}(V,F)$.
\end{thm}

\begin{pf} If for each $B\sim F$  and ${\cal L}_o$ transversal to $B$, we consider the representable
functors
$\fu{F_B}(S)= {Hom}_{{\cal O}_S}(({\cal L}_o)_S,\w B_S)$, then we obtain an open covering of
schemes for
$\fu{\text{Gr}}(V,F)$.
 \end{pf}

Fixing $V^+\sim F$, the index of an element ${\cal L}\in \fu{\text{Gr}}(V,F)(S)$, is defined as the
locally constant function
$i_{{\cal L}}\colon S\to\Bbb Z$ defined by:
$$i_{{\cal L}}(s)=\text{dim}_{k(s)}({\cal L}_s\cap \w {V^+}_s)-\text{dim}_{k(s)} \w V_s/{\cal L}_s+ \w
{V^+}_s$$

\begin{rem}{\label{ind}}
Let $\text{Gr}^d(V,F)$ be the moduli scheme of discrete submodules  of index $d$. One then has: 
$$\text{Gr}(V,F) =\underset{d\in\Bbb Z} \coprod\text{Gr}^d(V,F)$$
\end{rem} 
\bigskip
For further  details about this construction  see  \cite{AMP}.

For our purpose $F'=\{I(\w O^\infty)^n\}_{I\subset A}$ and $V=(\w{\Bbb A}^\infty)^n$ ($I$ ideal of
$A$).
We denote now  $\text{Gr}((\w{\Bbb A}^\infty)^n,F')$   by
$\text{Gr}((\w{\Bbb A}^\infty)^n,(\w O^\infty)^n)$

 
\section{The  Krichever morphism as   a closed immersion   in the Sato's infinite Grassmannian.
\label{grass}}

 Let $S$ be   scheme over $k$ ($k$ not necessarily finite), $X$ a proper, smooth and geometrically
irreducible over $k$, 
$\pi:X\times S
\to S$  the natural projection, and $M$ a  locally free sheaf on $X\times S$ of rank
$n$.  We shall denote by $O_M^\infty$   the sheaf on $S$, $\pi_*(i_*(M_{\vert
Spec({\cal O}_\infty)\times S}))$, where ${\cal O}_\infty$ is the stalk of ${\cal O}_X$ in
$\infty \in X$ and $i:Spec({\cal O}_\infty)\times S\to X\times S$ is the natural
morphism. For an effective divisor    $D$   
  on $X$, such that $\infty \notin supp(D)$, ${\cal O}_{X\times S}(D)$ will be denote by 
${\cal O}_X(D)
\underset k
 {\otimes}{{\cal O}}_S $.

\begin{defn}A $D$-level structure on $ S$, $( M,f_D)$, is a surjective
morphism of sheaves of ${\cal O}_{X\times S}$-modules

$$f_D:M\to    ({{\cal O}_X}/{{\cal O}_X(-D)})^n\underset k  {\otimes}{{\cal O}}_S $$  
 \end{defn}

Two $D$-level structures, $(M,f_D)$ and $(M',{f_{D'}} )$, are said to be equivalent  if there
exists an isomorphism of sheaves of ${\cal O}_{X\times S}$-modules, $\phi:M\to M'$, such that the
following diagram is commutative.

$$  \CD M @>\phi>>
 M' \\ @ V{f_D}VV    @ V{{f'}_D}VV \\ (O_D^n)_S @>=>Id>  (O_D^n)_S
\endCD$$   

being $(O_D^n)_S= ({{\cal O}_X}/{{\cal O}_X(-D)})^n\underset k  {\otimes}{{\cal O}}_S $

Let us call $\fu{\cal M}_D^n$  the moduli functor of  $D$-level structures.\bigskip

\begin{defn} A $\infty'$-formal level structure on $S$, $(M,f^\infty)$, is an
element of 
 
$$\w{\fu{ \cal M}}_\infty^n(S)=\underset
{\infty \notin supp(D)}   {\underset {D>0}
 \varprojlim}\fu {\cal M}^n_D(S)$$
\end{defn}

We denote by $\w{{\cal O}}_S^\infty $ and $\w{\cal A}_S^\infty $  the sheaves of  ${\cal O}_{X\times
S}$-modules:  
$$\underset {\infty \notin supp(D)} {\underset {D>0} \varprojlim} 
\mod {  {\cal O}_X}/{{\cal O}_X(-D)}\underset k {\otimes}{{\cal O}}_S \quad and \quad   
 \w{{\cal O}}_S^\infty \underset {\o_{X\times S}}  \otimes i_*({\cal O}_\infty \otimes {\cal O}_S)
 $$
respectively.  

If we take the direct image   by $\pi:X\times S \to S$ in $\w{{\cal O}}_S^\infty $ and 
$\w{\cal A}_S^\infty$, we obtain two sheaves of $\o_S$-modules, which will be called $\w{ 
O}_S^\infty $ and $\w{ 
\Bbb A}_S^\infty $. In the case $S=Spec(k)$,  $\w{{\cal O}}^\infty $    is
precisely  $\w{  O} ^\infty $ and $\w{\cal A}^\infty$ is $\w{\Bbb A}^\infty $.
 
If $(M,f^\infty)$  is a $\infty'$-formal level structure on $S$ we can obtain an exact secuence
of sheaves of ${\cal O}_{X\times S}$-modules:

$$0\to M
 \longrightarrow  i_* ( M_{\vert Spec({\cal O}_\infty) \times S} )  \overset{ \bar{f^\infty}} 
\longrightarrow   
     (\w{\cal A}_S^\infty    /     \w{{\cal O}}_S^\infty)^n
\to 0 $$

being  $ {f^\infty}$ the morphism obtained from the level structure and $\bar{f^\infty}$    the
composition of the morphisms: 
$$   f^\infty \otimes 1:M\underset {{\cal O}_{X\times S}}  \otimes  i_*({{\cal O}_\infty \otimes {\cal O}_S})=i_* ({M_{\vert Spec({\cal O}_\infty) \times S}})\longrightarrow (\w {{\cal O}}_S^\infty)^n\underset
{{\cal O}_{X\times S}}  \otimes  i_*({{\cal O}_\infty \otimes {\cal O}_S})= (\w{\cal A}^\infty_S)^n$$
and the natural projection  $ (\w{\cal A}_S^\infty )^n\to   {(\w{\cal A}_S^\infty)^n} /
{   (\w{{\cal O}}_S^\infty)^n  }$.

It is easily obtained   tensoring by $M$ the sequence of   ${\cal O}_{X\times
S}$-modules: 
$$ \quad 0\to
{\cal O}_{X\times S}
\overset {h}    \longrightarrow   i_*({\cal O}_\infty \otimes {\cal O}_S)
\overset {j}     \longrightarrow   
 \w{\cal A}_S^\infty/ \w{{\cal O}}_S^\infty\to 0 $$
 where $h$ and $j$ are the natural morphisms.

So by tensoring by ${\cal O}_X(D)$ this exact sequence
and  taking $\pi_*$,   we obtain: 
 
$$\pi_*M( D)= O_\infty^M \cap ( \w{  O}_S^\infty (D))^n \quad\quad
   \quad R^1\pi_*M( D)=\mod{ ( \w{\Bbb A}_S^\infty) ^n   } / { O_\infty^M+ ( { \w{  O}_S^\infty(
D)} )^n} $$
   $\w{  O}_S^\infty (D)$ being $\pi_*(\w{ {\cal O}}_S^\infty \underset {{\cal O}_{X\times S}} 
\otimes
{\cal O}_{X\times S}(D))$.

\begin{lem}\label{kri} There exists a natural morphism of functors in the category of
$k$-schemes:
$$\varphi: \w{\fu{ \cal M }}_\infty^n 
\to 
\fu {\text{Gr}}((\w {\Bbb A}^\infty)^n,(\w {O}^\infty)^n)$$
being   $\varphi(M,f^\infty)  = \pi_*(f^\infty\otimes 1)(O_M^\infty)$. 
\end{lem}
\begin{pf}This morphism can be obtained  adapting the proof of \cite{M}, \cite{SW} using the
following vanishing theorem.
\end{pf}

\begin{thm}If $M$ is a  locally free sheaf  of rank $n$ on
$X\times S$, there exists    an affine  open covering $\{U_i\}_{ i\in I}$ of $S$ and   effective
divisors $D_i$ over $X-\{\infty\}$, with $deg({D}_i)>>0$, such that: $R^1\pi_*M(D_i)_{\vert U_i}=0$
and  
$\pi_*M(D_i)_{\vert U_i} $ is a free sheaf of ${\cal O}_{U_i}$-modules of finite type. Moreover, there
exists $-D'_i  \subset X-{\infty}$ with $deg(D'_i)<<0$ such that $\pi_*M(D'_i)_{\vert U_i}=0 $
\end{thm}
\begin{pf} Classically this result is proved by assuming that $S$ is noetherian. Since   
  an exact sequence of modules over a ring 
$$0\to M'\to M \to M''\to 0$$
 is split when $M''$ is free and if  $M$ is
free $M'$ is locally free, it is possible to avoid the noetherian hipothesis. 
\end{pf}

 \bigskip

\begin{thm}  Let $\cal F$ be the subfunctor of
the infinite Grassmannian  defined by:
$${\cal F}(S)= \{{\cal L}\in
\fu{\text {Gr}}(({\w{\Bbb A}_S^{\infty}})^{n}, (\w{ O}^{\infty} _S  )^{n})(S)\quad:\quad
\text{are } {\cal O}_\infty \text{-modules}\}$$ 

There then exists a natural morphism of functors: $T:{\cal F}\longrightarrow
 \w{\fu{ \cal M}}_\infty^n$, such that $T\circ \varphi
=Id_{\w{\fu{  \cal M}}_\infty^n}$ and $\varphi\circ T=Id_{\cal F}$.
The definition of $T$ will be clear from the proof.
\end{thm}

\begin{pf} If ${\cal L}$ is a ${\cal O}_\infty\otimes {\cal O}_S$-module,  ${\cal L}$ has an associated  
quasicoherent ${\cal O}_{X\times S}$-module $i_*\tilde {{\cal L}}$, where $\tilde{{\cal L}}$ is the
quasicoherent module over $Spec({\cal O}_{\infty})\times S$ given by ${\cal L}$.  
From the definition of $\tilde{{\cal L}}$ it is possible to deduce   a morphism of sheaves of
${\cal O}_{X\times S}$-modules:
$$g:i_*\tilde {{\cal L}}\to    (\w{\cal A}_S^\infty)^n$$

On the other hand, since   ${\cal L}$ is a discrete submodule, there exists  an affine open covering
$\{U_i\}_{i\in I}$ of $S$, and $D_i$ effective divisors over $X-\{\infty\}$ such that  
$$   (\w{\cal A}_{U_i}^\infty)^n /   (\w{{\cal O}}_{U_i}^\infty(D_i))^n+ i_*\tilde {{\cal L}}_{X\times U_i}  =0$$ 

Note that ${{\cal O}_{X\times S}}(D_i)$ are faithfully flat ${\cal O}_{X\times S}$-modules, and hence 
$$\bar g:i_*\tilde {{\cal L}}\to   {(\w{\cal A}_S^\infty)}^n / {   (\w{{\cal O}}_S^\infty)^n  }$$
 is a surjective morphism of quasicoherent sheaves of ${\cal O}_{X\times S}$-modules.
  
Let $M$ be the quasicoherent sheaf defined by $Ker(g)$. $M$ has   associated in a natural way a
$\infty$-formal level structure, $f^\infty$. This is deduced for each  $D\subset X-\{\infty\}$  
from the commutative diagrams:

$$\text{(3.5.1} \CD 0  @>>> M @>>> i_*\tilde {{\cal L}} @>\bar g>> (\w{\cal A}_S^\infty)^n   @>>>0
0  \\@. @AAA @AA=A @AAA @.  \\
  0 @>>>M(-D) @>>> i_*\tilde {{\cal L}} @>\bar g(-D)>>  (\w{{\cal O}}_S^\infty(-D))^n  
@>>> 0
\endCD$$
 and the snake Lemma.

To finish the proof, the two following Lemmas will show that $M$ is locally free of
rank
$n$ over
$X\times S$. Then, $T({\cal L})$ will be defined by $(M,f^\infty)$.

From the very definition of $T$ and $\varphi$, it is not hard to prove that $T\circ \varphi
=Id_{\fu{ \w{\cal M }}_\infty^n}$ and $\varphi\circ T=Id_{\cal F}
$.   \end{pf}\bigskip \bigskip

(3.5.2)  From the exact sequence 
$$ \CD 0  @>>> M @>>> i_*\tilde {{\cal L}} @>\bar g>> (\w{\cal
A}_S^\infty)^n   @>>>0 \endCD$$

    for $deg(D')>>0$, $\pi_*M(D')$ is locally free on $S$, $M(D')$ is locally
generated by global sections on $X\times S$  and
$R^1\pi_*M(D')=0$. Moreover,   $M_s=M\underset {{\cal O}_{X\times S} } \otimes({\cal O}_X\otimes
k(s))$, for each $s\in S$, is a  locally free sheaf of rank $n$    on $ {X\times k(s)}$ and $M$ is
flat over $S$, from this last staments 
 we deduce that if $S$ is locally  noetherian, $M$ is locally of
finite presentation and thus $M$ is locally free of rank $n$. And if $S$ is not locally noetherian,
we have to prove this by other  methods:

\begin{lem}$  M $ is a  locally free sheaf of rank $n$ on $ {X\times
S}$.\end{lem}

\begin{pf} We can assume that $k$ is algebraically closed, $S=Spec(B)$, $p\in  Spec(A)$
  verifying
$\quad H^1\left(X\times Spec(B),  M(-p)\right)  =0
$, $H^0\left(X\times Spec(B),  M(-p)\right)  $ is a $B$-module locally free of finite rank and that
there exists a surjective morphisms of sheaves 
$$(3.6.1){{\cal O}}_X \underset k \to \otimes H^0\left(X\times Spec(B),  M \right)    \to    M 
\to 0
$$  

Taking global sections in the exact sequence of the $p$-level structure on $M$ deduced from the
snake lema for (3.5.1) 
we obtain an exact sequence of $B$-modules 
 
$$(3.6.2)\quad 0\to H^0\left(X\times Spec(B),  M(-p)\right)   \longrightarrow
H^0\left(X\times Spec(B),  M \right)  
\overset {f_p}  \longrightarrow 
 B^n\to 0$$ 

Let $t_1,...,t_n$ be elements of $H^0\left(X\times Spec(B),  M \right)$ such that their
images by the morphism $f_p$ are a basis of $B^n$. Let us consider the natural morphism
    
$$ \quad t_1({\cal O}_X\otimes B) \oplus ...\oplus t_n({\cal O}_X\otimes B) \to M $$ 
This morphism is injective since that $M/m^k_pM=(O_p/m^k_pO_p)^n\otimes B$, for all $k\in {\Bbb
N}$. 

In conclusion, if for a given  $y\in X\times Spec(B)$ ($y$ is not necessarily of the form
$p\times s$) it suffices to find elements
$\quad t'_1,...,t'_n\in H^0\left(X\times Spec(B),  M \right)$ such that 

$$\left(\frac {M}{t'_1({\cal O}_X\otimes B) \oplus
...\oplus t'_n({\cal O}_X\otimes B)}\right)_{\vert{k(y)}}=0$$  
 (recall that $M$ is of finite type).  

By   (3.5.2), we obtain for each $y\in X\times Spec(B)$ the isomorphism $M\vert k(y)\simeq
k(y)^n$. Moreover, we have the surjective morphism (3.6.1) and by the exact sequence (3.6.2) 
 
$$H^0\left(X\times Spec(B),  M \right)
=(t_1B\oplus...\oplus t_nB)\oplus H^0\left(X\times Spec(B),  M(-p) \right)$$ 
Thus,   in the case of   $\{t_1,...,t_n\}$ not being a basis in 
  $M_{\vert {k(y)}}  $, as $k $ has enough elemnts, one could find    $\{z_1,...,z_n \}\in
H^0\left(X\times Spec(B),  M(-p)
\right)$ such that   $\{t_1+z_1,...,t_n+z_n\}$, would be a basis in  $M_{\vert {k(y)}}  $. Taking
$t'_i=t_i+z_i$, we conclude.\end{pf}  \bigskip \bigskip

We are now able to prove the main statement of this section, namely, that $\cal F$ (and therefore
$\w{\fu{ \cal M }}_\infty^n  $)
is a  closed  subfunctor of
$\text{Gr}((\w {\Bbb A}^\infty)^n,(\w {O}^\infty)^n)$.

\begin{thm} The functor $\w{\fu{ \cal M }}_\infty^n $ is representable by
a closed subscheme $ \w{ { \cal M }}_\infty^n  $ of

$\text{Gr}((\w {\Bbb A}^\infty)^n,(\w {O}^\infty)^n)$.\end{thm}
\begin{pf} The condition ${\cal L}'\subset {\cal L}$ defines a closed subscheme in $\text{Gr}((\w
{\Bbb A}^\infty)^n,(\w {O}^\infty)^n)$, taking $a.{\cal L}={\cal L}'$ for all $a\in {\cal O}_\infty$ we conclude.\end{pf} \bigskip \bigskip

 \begin{rm}  \rm{As consequence of Remark {\ref{ind}}, it is easy to check 

$$   \w{ { \cal M }}_\infty^n  =\underset {d\in \Bbb
Z}  \coprod    \w{ { \cal M }}_\infty^{n,d}   $$ 
where
$ \w{ { \cal M }}_\infty^{n,d}  $ is the moduli scheme
of vector bundles   of rank
$n$ and degree
$d+n(g-1)$,   with $\infty'$-formal
level structures.

Moreover, bearing in mind Lemma 3.4 of \cite{BL} and the above results, we obtain that 
$Sl^n(k((t)))/Sl^n(A)$ is a scheme in groups which is a closed subscheme of the infinite
Grassmannian
$\text{Gr}((k((t)))^n,(k[[t]])^n)$, being $t$   a local parameter of $\infty \in X$. In some
sense the ind-group scheme $Sl^n(A)$, is a parabolic subgroup of $Sl^n(k((t)))$.}\end{rm}


\section{Immersion of Drinfeld moduli schemes     in a finite product of
infinite  Grassmannians.} 
 Now $k$ is a finite field ${\Bbb F}_q$. Let $S$ be a ${\Bbb F}_q$ scheme.

\begin{defn}\label{drin}\rm{\cite{D1}, \cite{Mu}} An elliptic sheaf is a diagram of vector
bundles of rank 
$n$ over $X\times S$:

$$ \CD \dots  @>>> M_{-1} @>i_0>> M_{0} @>i_{1}>> \dots  @>i_n>>
M_{n } @>i_{n+1}>>     \dots \\@. @AA{t}A @AA{t}A @. @AA{t}A   \\
  \dots  @>>> F^* M_{-2} @>F^*i_{-1}>> F^*M_{-1} @>F^*i_{0}>> \dots  @>F^*i_{n-1}>>
F^*M_{n-1 } @>F^*i_{n }>>     \dots
\endCD$$

satisfiying:

a) For any $s\in S$, $deg((M_{-1})_s)=n(1-g)$.   $deg$ denotes the degree.  
 \bigskip \bigskip

b) For all $i\in \Bbb N$, $M_{i+n}=M_i(\infty)(=M\underset {{\cal O}_{X\times S}}   \otimes
({\cal O}_X(\infty)\underset k   \otimes {\cal O}_S))$. 
\bigskip \bigskip

c) $M_i+t(F^*M_i)=M_{i+1}$ ($t$ is a morphism of ${\cal O}_{X\times S}$-modules).

  $F$  denotes the morphism of schemes:
 $$ Id\times \sigma: Spec(A)\times S\longrightarrow Spec(A)\times S    $$ 
   $\sigma$ being the morphism of Frobenius on $S$. \end{defn}
 
  \bigskip

\begin{defn}\label{lest}   A $D$-level structure in a Drinfeld diagram 
  is a
$D$-level structure in each vector bundle $M_i$ compatible with the diagram.\end{defn}
It is not hard to see that if a Drinfeld diagram has a $D$-level structure, then for all
$i\in \Bbb Z$ is 
$\left(M_i/t(F^*M_{i-1})\right)_{\vert  D\times S}=0$

\begin{rem}\label{dri}\rm{It is known that there exists an equivalence of categories, which
commutes with base change:

 $$   \left\{\vcenter{\hsize=5.5cm 
\parindent=0pt
Drinfeld modules of rank $n$ over $S$
\par whith characteristic away from  
\par $supp(D)$ with an $I_D$-level structure }\right\} \leftrightarrow 
\left\{\vcenter{\hsize=5.5cm 
\parindent=0pt
Elliptic sheaves  of rank $n$ over $S$
\par s.t. $\left(M_i/t(F^*M_{i-1})\right)_{\vert   D \times S}=0$
\par   for all $i$, with a twisted $ D$-level structure }\right\}$$}\end{rem}

A twisted $D$-level structure  on a locally free $\o_{X\times S}$-module $M$, is a surjective
morphism of $\o_{X\times S}$-modules

$$f_D:M\to (\Omega^1_X/\Omega^1_X(-D))^n\underset k \otimes \o_S$$
( $\Omega^1_X$ is the sheaf of differentials of $X$).

For comodity in the notation,   we are going to work with the definition of level structure given
in
\ref{lest}, with a easy modification can be changed this results for twisted $D$-level structures.

 Regarding the Remark \ref{dri} and the results of the above sections, we shall show that the
moduli functor of Drinfeld $A$-modules with $\infty'$-formal level structures  is
representable by a subscheme  of a finite product of infinite Grassmannians.

Let us consider the moduli functor:
$$ \w{\fu{ \cal D }}_\infty^n =\underset {\infty \notin supp(D)} {\underset {D>0 } \varprojlim}
{\fu{
\cal D}^n_D } $$
where $ {\fu{ \cal D}^n_D } $ is the moduli functor of elliptic sheaves 
  with  level structures on the effective divisor $D$.

We want to prove:

\begin{thm} There exists an injective   morphism of functors in the category
of ${{\Bbb F}_q}$-schemes:

$$  \psi:\w{\fu{ \cal D }}_\infty^n \hookrightarrow \underset {i=0}   {\overset   { n-1}  
\prod} 
 \fu{\text{Gr}}^{i+1}((\w {\Bbb A}^\infty)^n,(\w {O}^\infty)^n)$$
such that $\w{\fu{ \cal D }}_\infty^n   $ is a locally closed subfunctor of $\underset {i=0}  
{\overset   { n-1}  \prod} 
 \fu{\text{Gr}}^{i+1}((\w {\Bbb A}^\infty)^n,(\w {O}^\infty)^n)$. So $\w{\fu{ \cal
D}}_\infty^n  $ is representable by a scheme $\w{  \cal D}_\infty^n$ of the infinite Grassmannian
$\underset {i=0}
  {\overset   { n-1}  \prod} 
  {\text{Gr}}^{i+1}((\w {\Bbb A}^\infty)^n,(\w {O}^\infty)^n)$.\end{thm}
\begin{pf}
To prove this Theorem we need a previous result about $\infty'$-formal level structures:

   If $(M,f^\infty)$ and $(M',{f'}^\infty)$ are two $\infty$-formal
level structures on $X\times S$, ${\cal L}$ and ${\cal L}'$ the discrete submodules associated with this
$\infty$-formal level structures and

$g: (M,f^\infty)\longrightarrow (M',{f'}^\infty)$ a morphism  between this level structures 

 i.e.: a commutative diagram of sheaves of ${\cal O}_{X\times S}$-modules

$$  \CD M @>g>>
 M' \\ @ V{f^\infty}VV    @ V{{{f'}^\infty}_D}VV \\ (\o_S^\infty)^n @>=>Id>  (\o_S^\infty)^n
\endCD$$

Then $g$ is unique. Moreover, the morphism $g$   exists if and only if ${\cal L}\subseteq {\cal L}'$. 

Now we are going to prove the theorem.  If $H$ denotes the Drinfeld diagram together with a
$\infty$-level structure:

$$ \CD \dots  @>>> (M_{0},f^\infty_{0}) @>>> \dots @>>> (M_{n },f^\infty_n) @>>>\dots 
\\@.@AA{t}A @.  @AA{t}A @.  \\
  \dots  @>>> (F^*M_{-1},F^*f^\infty_{-1}) @>>> \dots @>>> (M_{n-1 },f^\infty_{n-1}) @>>>\dots
\endCD$$

the morphism $\psi$ is defined as follows
$$\psi(H)=\left((\varphi(M_0,f^\infty_0),\cdots,(\varphi(M_{n-1},f^\infty_{n-1})\right)$$
( $\varphi$ being the Krichever morphism (\ref{kri})).

The injectivity of $\psi$ follows easily from b) of (\ref{drin}) and noting that  $g$
is fixed for the
$\infty'$-formal level structure.

On the other hand, the image functor of $\psi$  lies in $ \underset
{i=0}  {\overset   {n-1}
 \prod}
\w{\fu{ \cal M }}_\infty^{n,i+1}$, 
and by  the above statement and  (\ref{drin}), the necessary and sufficient condition for   an
element
$$({\cal L}_1,\cdots,{\cal L}_n)\in 
 \underset {i=0}  {\overset   {n-1}  \prod} 
 \w{\fu{ \cal M }}_\infty^{n,i+1} (S) \subseteq
\underset {i=0}  {\overset   {n-1}  \prod} 
 \fu{\text{Gr}}^{i+1}((\w {\Bbb A}^\infty)^n,(\w {O}^\infty)^n)(S)$$ 
to be  in $Img(\varphi)(S)$ is that:

\begin{itemize}
 
\item  For all $1\leq i\leq n-2$, ${\cal L}_i\subset {\cal L}_{i+1}$  and 
${\cal L}_{n-1}\subset {\cal L}_0(\infty)$, ${\cal L}_0(\infty)$ being 
  the discrete submodule obtained from the
$\infty$-formal level structure:

 $ (M_0\underset {{\cal O}_{X\times S}}  \otimes
{\cal O}_{X\times S}(\infty),f^\infty_0\otimes 1)$.

 \item  For all $1\leq i\leq n-2$, $(\sigma)^*{\cal L}_i\subseteq {\cal L}_{i+1}$ and $(\sigma)^*{\cal L}_{n-1}\subset {\cal L}_0(\infty)$.

\item For all $1\leq i\leq n-2$,  $(\sigma)^*{\cal L}_i+{\cal L}_i={\cal L}_{i+1}$ and $(\sigma)^*{\cal L}_{n-1}+{\cal L}_{n-1}={\cal L}_0(\infty)$

\end{itemize}

 Recalling that given two discrete submodules ${\cal L}$ and ${\cal L}'$ in $(\w {\Bbb
A}_S^\infty)^n$ the subset of $S$ where ${\cal L}\subseteq {\cal L}'$ is a closed subscheme of $S$,
then an element of 
$
\underset {i=0}   {\overset   {n-1}  \prod} 
 \w{\fu{ \cal M }}_\infty^{n,i+1} (S)$  fulfills
conditions 1) and 2) is a closed subscheme
 and by  Nakayama's Lemma, applied to ${\cal L}_{i+1}/{\cal L}_i+\sigma^*{\cal L}_i$,  
condition 3)
is fullfiled in an open subscheme of
$S$.

( ${\cal L}_{i+1}/{\cal L}_i$ is a finite $B$-module since  that there exists  a   surjective
morphism of
$B$-modules $ {\cal L}_{i+1}/{\cal L}_{i+1}(-\infty)\to {\cal L}_{i+1}/{\cal L}_i$. Recall that ${\cal L}_{i+1}/{\cal L}_{i+1}(-\infty)\simeq\pi_*(M_{i+1}/M_{i+1}(-\infty))$ is a $B$-module locally free
of finite rank). We thus conclude.\end{pf}


\vskip1.5truecm { \'Alvarez V\'azquez, Arturo}\newline {\it
e-mail: } aalvarez@@gugu.usal.es
\end{document}